\documentclass[aip, amsmath,amssymb, reprint, ]{revtex4-1}
\usepackage[utf8]{inputenc}
\usepackage[T1]{fontenc}
\usepackage[english]{babel}
\usepackage{amsmath}
\usepackage{graphicx}		
\usepackage{natbib}
\usepackage{textcomp}
\usepackage{gensymb}
\usepackage[hidelinks]{hyperref}
\usepackage{xcolor}
\usepackage{siunitx}
\usepackage{mathrsfs}
\usepackage{multirow}
\usepackage{amsthm}
\usepackage{float}
\usepackage{textgreek}
\theoremstyle{definition}

\graphicspath{{figs/}}
\usepackage{dcolumn}
\usepackage{bm}

\usepackage{mathptmx}
\usepackage{etoolbox}

\makeatletter
\def\@email#1#2{%
 \endgroup
 \patchcmd{\titleblock@produce}
  {\frontmatter@RRAPformat}
  {\frontmatter@RRAPformat{\produce@RRAP{*#1\href{mailto:#2}{#2}}}\frontmatter@RRAPformat}
  {}{}
}%
\makeatother

\begin{document}
\preprint{AIP/123-QED}

\title{Transients versus network interactions give rise to multistability through trapping mechanism}

\author{Kalel L. Rossi}
\affiliation{Theoretical Physics/Complex Systems, ICBM, Carl von Ossietzky University Oldenburg, Oldenburg, Lower Saxony, Germany}

\author{Everton S. Medeiros}
 \affiliation{São Paulo State University (UNESP), Institute of Geosciences and Exact Sciences, Rio Claro, SP, Brazil}

\author{Peter Ashwin}
\affiliation{Department of Mathematics and Statistics, University of Exeter, Exeter EX4 4QF, United Kingdom}

\author{Ulrike Feudel}
\affiliation{Theoretical Physics/Complex Systems, ICBM, Carl von Ossietzky University Oldenburg, Oldenburg, Lower Saxony, Germany}

\begin{abstract}
In networked systems, the interplay between the dynamics of individual subsystems and their network interactions has been found to generate multistability in various contexts. Despite its ubiquity, the specific mechanisms and ingredients that give rise to multistability from such interplay remain poorly understood. In a network of coupled excitable units, we show that this interplay generating multistability occurs through a \textit{competition} between the units' transient dynamics and their coupling. Specifically, the diffusive coupling between the units manages to \textit{reinject} them in the excitability region of their individual state space and effectively trap them there. We show that this trapping mechanism leads to the \textit{coexistence} of multiple types of oscillations: periodic, quasiperiodic, and even chaotic, although the units separately do not oscillate. Interestingly, we show that the attractors emerge through different types of bifurcations - in particular, the periodic attractors emerge through either saddle-node of limit cycles bifurcations or homoclinic bifurcations - but in all cases the reinjection mechanism is present. 
\end{abstract}

\maketitle

\begin{quotation}
A common behavior in nonlinear dynamical systems is multistability, the coexistence of multiple stable solutions. Multistability is observed both in nature and models in a wide variety of applications, including the climate, power grids, ecology, and the brain \cite{feudel2018multistability}. It has important consequences: a multistable system operating on a particularly desirable attractor may not be safe, as a perturbation in the state of the system can cause it to switch to another coexisting attractor. On the other hand, coexistence of attractors may be useful, for instance by enabling the implementation of memory functions \cite{khona2022attractor}. In networked systems, multistability can arise from the interactions of the multiple subunits, but the specific mechanisms that cause this emergence are still not fully known. In this work we demonstrate one mechanism that gives rise to multistability from the interaction of even only two units. The particular system we study is composed of excitable neurons coupled diffusively, as a simple model for neurons coupled through gap junctions. Importantly, the neurons individually do not oscillate, but the coupling between them leads to a rich variety of attractors with oscillations. We show in a network of ten neurons the emergence of periodic, quasiperiodic and even chaotic stable oscillations, which all emerge from the same underlying mechanism.
\end{quotation}

\section{Introduction}

The long-term behavior of dynamical systems is determined by their attractors, which are stable states that attract certain sets of initial conditions. Dynamical systems can possess several attractors coexisting for the same parameters, such that different initial conditions can lead to different long-term behaviors - a phenomenon called \textit{multistability} \cite{pisarchik2014control, feudel2018multistability}. In power grids, this can mean the difference between the proper functioning of the grid and a blackout \cite{motter2013spontaneous}; in ecological systems, it can mean the difference between extinction of a certain species and their survival \cite{meng2022the}. In neuronal circuits, multistability has been shown to be important for computations \cite{driscoll2024flexible}, and may, for instance, implement memory storage if the attractors correspond to different memories \cite{foss1996multistability, wilson1972excitatory}.

Many systems display multistability, particularly networked systems, in which individual units are coupled together according to some type of interaction \cite{pisarchik2014control}. An important type of interaction in networked system is diffusion. One example is found in interacting ecological patches, in which each patch has its own dynamics but also interacts with other patches by migration, or diffusion, of species \cite{pilosof2017the}. Another example is found in neuronal networks, in which neurons interact with each other through the transport of ions across their cell membrane \cite{sohl2005expression, benett2004electrical}. In these two examples, the interaction between units $i$ and $j$ can be modeled by a linear diffusion term dependent on the difference $x_j - x_i$ between the state variables $x_i$ and $x_j$ of units $i$ and $j$ \cite{stankovski2017coupling, loppini2015mathematical,kepler1990the, blasius1999complex, medeiros2021asymmetry, liang2023a, sadykov2021model, kaneko1992overview}. Understanding the emergence of multistability in networked systems with this kind of interaction therefore finds applications in many fields, and is still an area of active research. 

There is a wide literature studying networked systems whose units separately oscillate and which are diffusively coupled \cite{winfree1980the, winfree1967biological, berner2021generalized, pikovsky2001synchronization, rodrigues2016the, mehrabbeik2023synchronization, kaneko1992overview, kaneko1985spatiotemporal, kaneko1992global}. Multistability in these systems is also well-known, with the emergence of different types of coexisting attractors \cite{ullner2007multistability, ullner2008multistability, rossi2022shifts, yanagita2005pair, koseska2013oscillation, crowley1989experimental, sporns1987chaotic, aronson1987an, mondal2021spatiotemporal, ansmann2016selfinduced, kaneko1992global, kaneko1997dominance, kaneko1993chaotic}. For instance, Ref.~\onlinecite{ullner2008multistability} studied two coupled repressilators, 7-dimensional units that have stable oscillations when uncoupled, and find the emergence of different types of attractors. Some attractors have two units oscillating with a large amplitude and some have one unit at a large amplitude and another with a very small amplitude, called inhomogeneous limit cycles. When more units are coupled in a big network with $N=100$ units, the authors showed in Ref.~\onlinecite{ullner2007multistability} that a large number of such attractors can coexist. In coupled mechanical oscillators, two coupled rotors have also been shown to exhibit large multistability (more than 3000 attractors) \cite{feudel1998dynamical}.

Less is known about multistability when the units individually do not oscillate, although it is known that oscillations can still arise due to the coupling. An important work in this direction is due to Smale in 1976 based on an idea by Turing in 1952 \cite{turing1952the, Smale1976, stankovski2017coupling}. Smale proposed the emergence of oscillations from non-oscillating units which have only one equilibrium that is stable and globally attracting in a region of their state space. It was shown that the oscillations come from a Hopf bifurcation, in which the equilibrium becomes unstable and a stable oscillation emerges \cite{Smale1976, pogromsky1999on}. Important work has also been done by Winfree, showing the emergence of waves in continua of identical excitable systems coupled through nearest-neighbor diffusive interactions \cite{winfree1980the}. Chaotic oscillations can also emerge from diffusive coupling applied to units with a single stable equilibrium in a region of state space. An example was given in Ref.~\onlinecite{kocarev1995on} for two coupled Chua circuits. Recently, researchers provided rigorous conditions for the emergence of chaos due to diffusive coupling \cite{nijholt2023chaotic}. However, these works generally do not look at multistability. Furthermore, they deal with a single equilibrium in a region of state space, and have not yet looked at a scenario in which more invariant sets, such as unstable equilibria, may also play a role.

The presence of unstable equilibria can alter the transient dynamics of non-oscillating systems. In some classes of models, which we study here, the unstable equilibria lead to a type of excitability \cite{izhikevichbook}. In this case, the unstable equilibria force part of the trajectories to go through a long excursion in state space, called an \textit{excitation}, before reaching the stable equilibrium. These excitations are common in neuronal models, where they correspond to a neuronal spike \cite{izhikevichbook}. Reference \onlinecite{yanagita2005pair} has described multistability emerging in two excitable FitzHugh-Nagumo neurons that were coupled repulsively, but over a relatively small parameter range. For attractive coupling, the authors did not observe multistability. A similar scenario was reported in Ref.~\onlinecite{ronge2021splay}, which studies excitable phase oscillators near a saddle-node bifurcation on an invariant circle. They show that coupling these oscillators with a repulsive diffusive term can generate stable periodic solutions, in particular splay states and cluster states.

In this work, we present two findings. First, we show that an attractive diffusive coupling can indeed create new attractors in coupled excitable systems. In fact, a wide variety of them: periodic, quasiperiodic and even chaotic oscillations arise by coupling excitable units, with $N=2$ units already being sufficient for periodic and quasiperiodic attractors. For larger networks of $N=10$ units, we show that these attractors, periodic, quasiperiodic and chaotic attractors can \textit{coexist} in the same system, for a range of coupling strengths, in line with results of strong multistability in networks with many units \cite{rossi2022shifts, ullner2007multistability, gelbrecht2020monte}.

The second finding contributes to an understanding of one mechanism through which these attractors emerge. We study their geometry, looking at the interaction between the units' local dynamics, which creates the excitability, and the diffusive coupling term, which pulls the units toward each other. We show that the competition between these two terms can trap the units in a particular region of their state space where excitability occurs. Based on this mechanism, the previously transient excitable dynamics is now repeatedly activated, generating permanent oscillations. This occurs for all the attractors observed, which emerge under different bifurcation scenarios. It also extends to networks with more than two interacting units, suggesting a powerful mechanism for the coexistence of a multitude of attractors in networked systems. 
 
A similar idea of attractors emerging when units are trapped in transient regions of state space has been previously reported in the literature \cite{medeiros2018boundaries, medeiros2019state, medeiros2021the, contreras2023scale, medeiros2023transient}. The mechanisms underlying this trapping were different, and multistability had not been reported. In Refs.~\onlinecite{medeiros2018boundaries, medeiros2019state, medeiros2021the} the authors study units with chaotic saddles (unstable chaotic sets) in their uncoupled state space. They show that the diffusive coupling manages to trap the units in that region, generating a chaotic motion that numerically appears to be stable. Meanwhile, authors in Ref.~\onlinecite{contreras2023scale} have shown trapping in the vicinity of canard transitions. By diffusively coupling Fitzhugh-Nagumo units with relaxation oscillations near a singular Hopf bifurcation, they showed that weak coupling leads to a metastable regime with trajectories switching between being transiently trapped near an unstable equilibrium - generating small amplitude oscillations - and escaping this trapping to perform a large excursion (excitation) in state space - generating a large amplitude oscillation \cite{contreras2023scale}. In these reported mechanisms, only one additional attractor emerges due to the coupling. An interesting property of the behavior we study in this paper is that the trapping generates a plethora of coexisting attractors. Therefore, although the local dynamics of the units is relatively simple, it can create rich multistable dynamics. As such, we believe it serves as a simple yet powerful example of the more widespread phenomenon of multistability through trapping. 

The paper is organized as follows. We describe the model and algorithms in Sec. \ref{sec:methods}. Then, Sec. \ref{sec:results} introduces the rich multistability seen in a network of $N=10$ units, from which we reduce to $N=2$ units to better understand the mechanism giving rise to this multistability. In Sec. \ref{sec:discussions} we then discuss these findings in relation to each other and to preexisting literature.

\section{Methods}\label{sec:methods}
\subsection{Model}\label{sec:uncoupled}
In this work we study networks formed by coupling two-dimensional units with state variables $x$ and $y$ whose evolution we write as:
\begin{align}\label{eq:generic-network-equation}
    \dot{x}_i &= f_1(x_i, y_i) + \epsilon_1 h_i(\mathbf{x})\\
    \dot{y}_i &= f_2(x_i, y_i) + \epsilon_2 h_i(\mathbf{y}),
\end{align}
with $\mathbf{x} = (x_1, \ldots, x_N) \in \mathbb{R}^N$ and $\mathbf{y} = (y_1, \ldots, y_N) \in \mathbb{R}^N$ the state variables of the system. We refer to $\mathbf{f}_i = \mathbf{f}(x_i, y_i) = (f_1(x_i, y_i), f_2(x_i, y_i))$ as the local dynamics of unit $i$ and to $\mathbf{h}_i  = (h_i(\mathbf{x}), h_i(\mathbf{y}))$ as the coupling term of unit $i$, which allows it to receive influence from other units. The parameters $\epsilon_1$ and $\epsilon_2$ control the strength of the interactions, with $\epsilon_1 = \epsilon_2 = \epsilon$ unless stated otherwise. The interaction is specified by a diffusive coupling of the form:
\begin{align}
  h_i(\mathbf{z}) = \sum_{j \in \Omega_i} (z_j - z_i),  
\end{align}
where $\mathbf{z}$ is either $\mathbf{x}$ or $\mathbf{y}$, and $\Omega_i$ is the set containing the indices $j$ of units connected to unit $i$, also called the neighborhood of $i$. 

For the local dynamics, we choose a simple two-dimensional model for a spiking neuron following the Hodgkin-Huxley formalism, as written by Ref. \onlinecite{izhikevichbook}. The dynamics of this model is described by the following functions:
\begin{align}
\label{eq1:neuron-model}
    f_1(x_i, y_i) &= \big(I - g_L (x_i - E_L) - g_{Na} m_\infty(x_i) (x_i - E_{Na})  \\
    &-g_K y_i (x_i - E_K) \big)/C, \notag\\
\label{eq2:neuron-model}
    f_2(x_i, y_i) &= (n_\infty(x) - y_i) / \tau,
\end{align}
where the neuron membrane potential and conductance variable are represented by $x$ and $y$, respectively.  The activation functions $m_\infty(x_i)$ and $n_\infty(x_i)$ are given by:
\begin{align}
 m_\infty(x_i) &= \frac{1}{1 + \exp (m_h - x_i)/k_m) }, \\
 n_\infty(x_i) &= \frac{1}{1 + \exp ((n_h -x_i)/k_n) }.
\end{align}

The parameters used are $\tau = \SI{0.16}{\milli\second}$, $C = \SI{1}{\micro\farad/\centi\meter\squared}$, $E_L = -\SI{80}{\milli\volt}$,$g_L = \SI{8}{\milli\siemens/\centi\meter\squared}$, $E_\mathrm{Na} = \SI{60}{\milli\volt}$, $g_\mathrm{Na} = \SI{20}{\milli\siemens/\centi\meter\squared}$, $E_\mathrm{K} = -\SI{90}{\milli\volt}$, $g_\mathrm{K} = \SI{10}{\milli\siemens/\centi\meter\squared}$, $m_h = \SI{-20}{\milli\volt}$, $k_m = \SI{15}{\milli\volt}$, $n_h = \SI{-25}{\milli\volt}$,  $k_n = \SI{5}{\milli\volt}$ and $I = \SI{2.0}{\micro\ampere/\centi\meter\squared}$. The dynamics of this system is very similar to that of the Morris-Lecar model \cite{morris1981voltage, izhikevichbook}. A slight increase in the membrane voltage $x$ leads to a quick increase in the Sodium current, which is negative ($(x-E_\mathrm{Na}) < 0$) and acts to increase the voltage even further in a positive feedback that rapidly increases $x$, initiating the excitation (spike). At sufficiently high voltage, the Potassium current increases, being activated by the conductance variable $y$. This current is positive ($(x-E_\mathrm{K})>0$) and becomes sufficiently large that it overcomes the Sodium current and decreases the voltage back to baseline, terminating the excitation and returning to the stable equilibrium. For a more in-depth explanation of the model and a complete explanation of the parameters, we refer the reader to Ref.~\onlinecite{izhikevichbook}. For simplicity, from now on we refer to the parameters without their corresponding units. We remark that technically the coupling term $\epsilon_1 h_i(\mathbf{x})$ acts as a current and should be divided by the capacitance $C$. For simplicity of notation, we consider $C$ to be already included in $\epsilon_1$. This is also not a problem since in this parametrization $C=1$.

\begin{figure}
    \centering
    \includegraphics[width=\columnwidth]{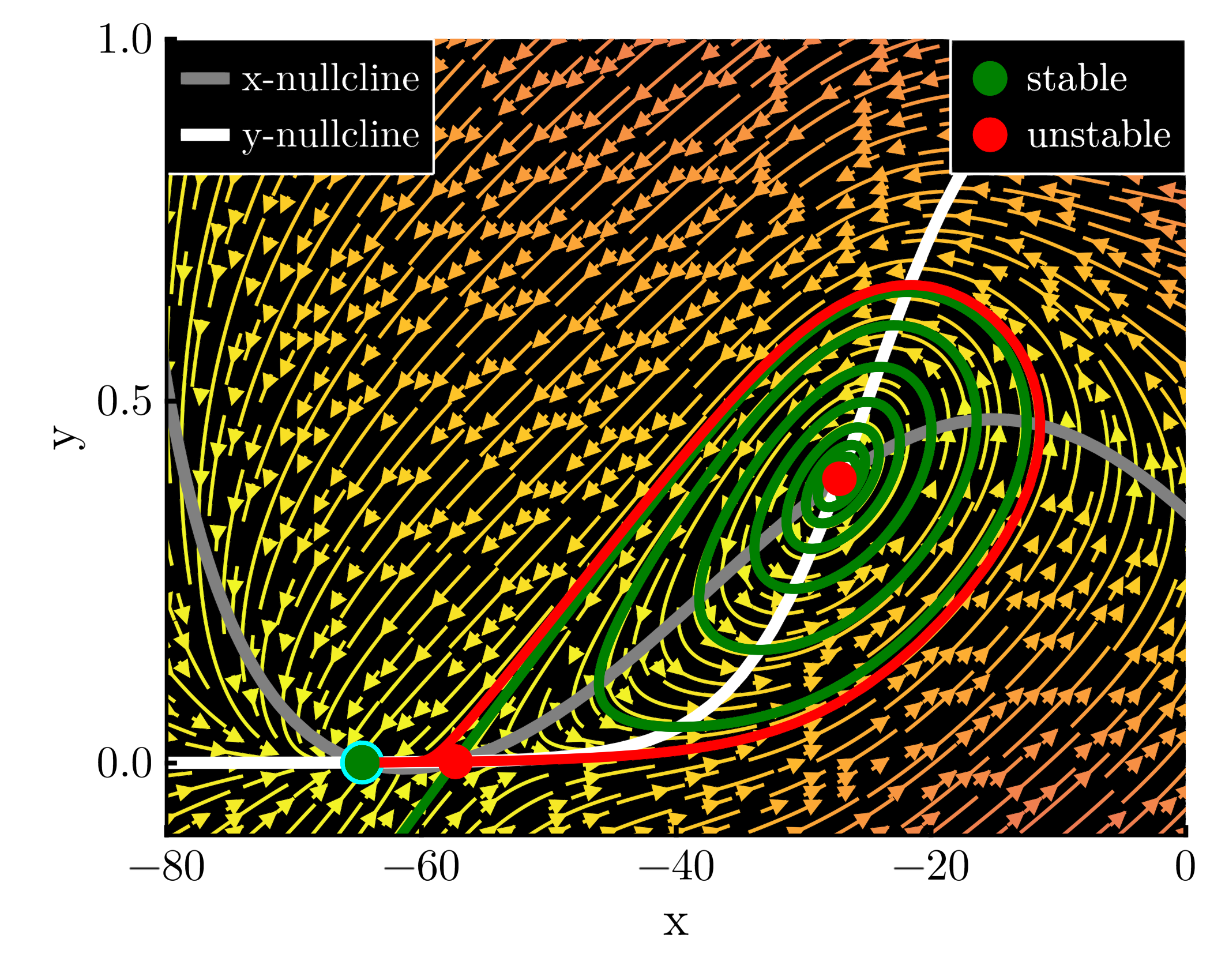}
    \caption{Phase portrait of the excitable uncoupled units. The green dot represents the stable node of the system at $(-63.3, 0.0005)$, the red dots represent the unstable focus at $(-27.1, 0.4)$ and the saddle point at $(-58.6, 0.001)$, with its stable and unstable manifold branches in the green and red lines. The $x$ and $y$ nullclines are indicated respectively in gray and white curves. The phase portrait is represented by the arrows, indicating the directions of the flow. As the flow indicates, there is a wide region in which trajectories must go around the stable manifold to reach the node. They correspond to a neuronal spike, since this is a sharp increase and then decrease in the membrane potential. This region is called the excitability region. Attractors emerging from the coupling live in this region of long transients. }
    \label{fig:excitable-neuron}
\end{figure}

For fixed $I = 2.0$, and the previously described parameters, the neuronal dynamics of the uncoupled units ($\epsilon_1 = \epsilon_2 = 0$) is excitable. The state space of the unit, shown in Fig.~\ref{fig:excitable-neuron}, is composed of a stable node (green circle), a saddle-point $\mathbf{x}_s^\mathrm{unc}$ (red circle close to the node), and an unstable focus (red circle). The stable manifold $\mathbf{W}^s(\mathbf{x}_s^\mathrm{unc})$ and the unstable manifold $\mathbf{W}^u(\mathbf{x}_s^\mathrm{unc})$ of the saddle are depicted in green and red lines, respectively. Additionally, the $x$-nullcline, defined as $\dot{x} = 0$, and the $y$-nullcline, defined as $\dot{y} = 0$, are shown in gray and white, respectively. 
As indicated by the vector field, the stable manifold $\mathbf{W}^s(\mathbf{x}_s^\mathrm{unc})$ roughly separates the state space into two regions: one that directly converges to the stable equilibrium, and another wherein trajectories go through long excursions before converging to the equilibrium. The long excursions are called excitations, and the region is called the \textit{excitability region}. 
%

In the main text of the manuscript, we focus on the phenomenology underlying the excitable case prescribed by $I=2.0$. In the Supplemental Material, we show that increasing $I$ leads to a homoclinic bifurcation, creating a stable limit cycle, followed by a saddle-node bifurcation that destroys the node and saddle of the units. We then discuss the effects of these bifurcations on the results presented in the paper.

\subsection{Numerical algorithms}

To find the attractors in our networked systems, we followed the method developed in Refs.~\onlinecite{datseris2023framework, gelbrecht2020monte, stender2021bstab}, which distinguishes between attractors based on user-defined features that uniquely characterize the attractors. To achieve this, it first integrates randomly chosen initial conditions in a specified region of the state space. The corresponding trajectories are then labeled based on their features, such as the mean value of their amplitude. These features must be chosen so that trajectories on different attractors exhibit distinct feature values. Subsequently, the features are separated using a grouping algorithm, which may involve clustering or simply distinguishing features that are more distant than a predefined threshold. This method works well for both low- and high-dimensional systems. We performed extensive numerical studies to ensure no attractors were missed, but this cannot be guaranteed. An attractor with a sufficiently small basin may, by chance, not be found. To reduce the risk of missing attractors, we used 5000 initial conditions for the $N=10$ results and 1000 for $N=2$. Test runs with more initial conditions did not find any further attractors. Each trajectory was integrated for a very long time, with a total transient time of 7000 and total integration time of 40000. The features used were the average pairwise Euclidean distance between the states of the units, their frequencies, amplitudes, and average position in state space. For equilibria, only the average position is considered, as the frequencies and amplitudes are zero.

The algorithms, with a complete documentation, are implemented in the Julia \cite{bezanson2017julia} package Attractors.jl \cite{datseris2022effortless, datseris2023framework}. We also verified the accuracy of results shown through continuation analysis using the XPPAUT 8.0 software \cite{ermentrout2002simulating}, finding the bifurcations giving rise to the attractors.

Integration was done with the package DifferentialEquations.jl \cite{rackauckas2016differential}, with the aid of packages DynamicalSystems.jl \cite{datseris2018dynamical} and DrWatson.jl \cite{datseris2020drwatson}. Plots were made with Makie.jl \cite{danisch2021makie}. The Tsitouras 5/4 Runge-Kutta method was used for the integrations, with absolute and relative tolerances of $10^{-9}$. The code for the analysis is publicly available in a GitHub repository \cite{rossi2024multistabilitygithub}.

\section{Results}\label{sec:results}

\subsection{Rich multistability of oscillations with 10 units}
\begin{figure*}[ht!]
    \centering 
    \includegraphics[width=\textwidth]{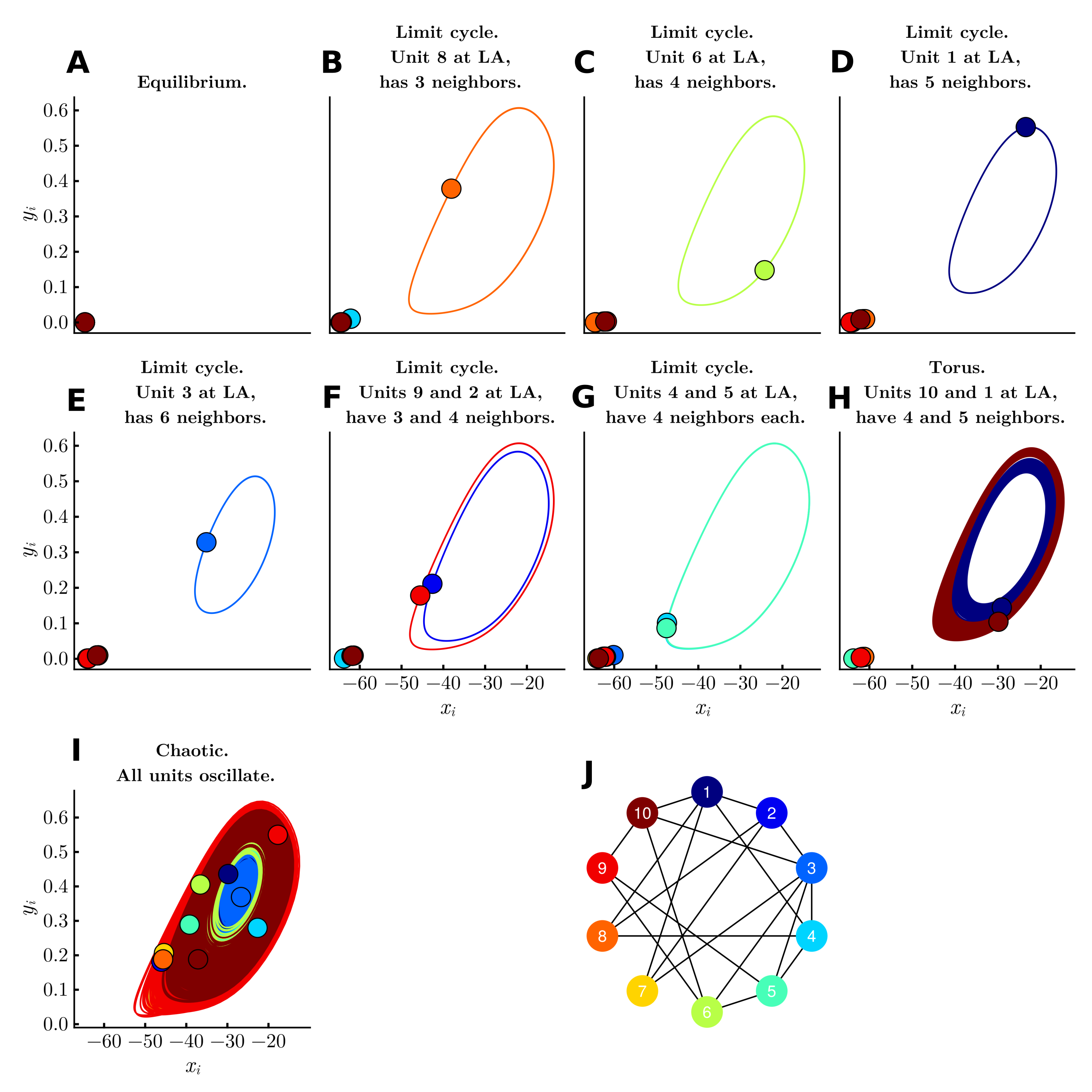}
    \caption{Rich multistability arising from diffusive coupling. Panels I-H show the stable equilibrium, periodic, quasi-periodic and chaotic attractors that coexist in the same network with $N=10$ randomly coupled units (shown in Panel J) with $\epsilon=0.15$. Each panel shows a trajectory on one of the attractors, projected onto the $x_i-y_i$ subspaces of each unit, all overlaid on top of each other. Circles correspond to the positions of the units at some arbitrarily chosen time point. The units (circles and curves) are colored from blue to red according to their index, as shown in the topology of panel J, such that unit $1$ is a deep blue and unit $10$ is a deep red.}
    \label{fig:multistability-network}
\end{figure*}

The diffusive coupling between the excitable units can generate rich oscillatory dynamics, in which equilibria coexist with periodic, quasiperiodic, and even chaotic oscillations. As we see later, these oscillations arise from the interplay between the diffusive coupling and the local flow field of the units. An example of these attractors is shown in Figs.~\ref{fig:multistability-network}A-I, in a network of $N=10$ excitable units following the topology shown in Fig.~\ref{fig:multistability-network}J. This topology is arbitrary and serves as an illustrative example, with the connections having been chosen at random. In Figs.~\ref{fig:multistability-network}A-I, we project the network state space into subspaces $x_i - y_i$ corresponding to each unit $i$, and overlay them all on top of each other.  In addition, the coupling strength is chosen as $\epsilon_1 = \epsilon_2 = \epsilon = 0.15$. 

The first type of attractor is shown in Fig.~\ref{fig:multistability-network}A. It corresponds to all units on the stable equilibrium, which is already present in the uncoupled units. This is the simplest solution, which must exist because, when the units are completely synchronized, the coupling term becomes zero and they follow their uncoupled dynamics, converging to the equilibrium.

The second type of attractor corresponds to one unit oscillating periodically with a large amplitude while the $N-1$ other units oscillate with a very small amplitude at a position between the stable equilibrium and the saddle of the uncoupled dynamics. The dynamics in this type of attractor resembles the so-called solitary states, since one unit behaves differently from the rest of the network. Such symmetry-broken solutions have been observed in regular \cite{semenova2018mechanism, jaros2018solitary, rybalova2019solitary,hellmann2020network}, adaptive \cite{berner2020solitary}, and complex networks \cite{schulen2022solitary}. For the chosen parameters, we have identified four stable solitary states, shown in Figs.~\ref{fig:multistability-network}B-E. The unit displaying a high-amplitude oscillation is said to be solitary. Interestingly, the amplitude of its oscillation is inversely proportional to the number of neighbors it has. With more neighbors, the coupling terms $\mathbf{h}_i(\mathbf{x})$ and $\mathbf{h}_i(\mathbf{y})$ of the solitary unit $i$ increase, and the amplitude of its oscillation decreases. The reason for this will become clearer in Section~\ref{sec:2units}, where we study in depth the case $N=2$. Numerical bifurcation analysis (not shown) reveals that these periodic attractors emerge in homoclinic bifurcations and disappear in saddle-node of limit cycle bifurcations (SNLC).

The third type of attractor corresponds again to periodic oscillations, but with two high-amplitude units, shown in Figs.~\ref{fig:multistability-network}F-G. In Fig.~\ref{fig:multistability-network}F, the units exhibiting high-amplitude oscillations are 2 and 9. Note that unit 2 has one more neighbor than unit 9, and its amplitude is smaller. In Fig.~\ref{fig:multistability-network}G, the units are 4 and 5. They have the same number of neighbors, so their amplitudes are identical. This type of attractor is thus a two-unit cluster periodic state. Bifurcation analysis reveals that these attractors emerge and disappear through SNLC bifurcations.

A fourth type of attractor also involves two units (1 and 10) oscillating with large amplitude, but now quasi-periodically, as shown in Fig.~\ref{fig:multistability-network}H. Similarly to the previous cases, the amplitude of their oscillations is proportional to the number of neighbors they have. As shown in the topology in Fig. \ref{fig:multistability-network}J, the oscillating units (1 and 10) are connected, so they also pull each other in directions perpendicular to their oscillations as they oscillate. Intuitively speaking, we can imagine that this interaction enlarges the width of the torus. Indeed, if one introduces a coupling parameter directly between units 1 and 10, i.e., setting $h_i(\mathbf{z}) = \sum_{\Omega_i}\epsilon_{i,j}(z_j-z_i)$, and specifically increasing $\epsilon_{1,10} = \epsilon_{10,1}$ from $0$ to $\epsilon$, the width of the torus in the $x_i-y_i$ projection increases. Thus, when oscillating together with different amplitudes, the coupling between the units causes their quasi-periodic curves to become broader. We see the emergence of tori in greater depth when we study the $N = 2$ case in Section~\ref{sec:2units}.

Finally, the fifth type of attractor involves all units oscillating together chaotically. All neuronal units are thus spiking chaotically in this attractor in a desynchronized fashion. The chaotic behavior, along with the periodic and quasiperiodic examples from earlier, has been verified by calculating the Lyapunov exponents of these attractors.

The results conveyed in Fig.~\ref{fig:multistability-network} occur for an intermediate range of coupling strength values, at and around $\epsilon=0.15$. Bigger coupling strengths tend to generate fewer attractors, and ultimately for strong coupling only the stable equilibrium exists. For weaker coupling, even more attractors can appear. In fact, for a range roughly between $\epsilon=0.05$ and $\epsilon=0.1$, more than 50 attractors can be found. These correspond to the various combinations of units having a very small amplitude oscillation, and units having a large amplitude oscillation. 

Furthermore, in the network we analyzed so far the units are coupled in both $x$ and $y$ directions (i.e., $\epsilon_1=\epsilon_2=\epsilon$). If only the $x$-direction is coupled ($\epsilon_1=\epsilon$, $\epsilon_2=0$), there still is multistability, but with fewer attractors. The $x$-coupling tends to stabilize attractors with more units oscillating at a higher amplitude, such that one can have 5 units oscillating periodically at a large amplitude and 5 with small amplitude, for instance.

To summarize, the addition of a simple linear interaction through the attractive diffusive coupling creates a plethora of oscillations from non-oscillating units in an excitable regime. The coupling is clearly able to counteract the units' tendency to converge onto the stable equilibrium. Our goal in the following sections is to elucidate this mechanism in more detail. To achieve this, we simplify our system and reduce the problem to $N=2$ interacting units.

\subsection{Emergence of attractors in a two-unit network} \label{sec:2units}
To illustrate the effect of the diffusive coupling on the excitable neurons, we show the attractors of the system for different coupling strengths for $N=2$ coupled units. Similarly to Fig.~\ref{fig:multistability-network}, each panel in Fig.~\ref{fig:attractors-2-units} shows the variables $x_i-y_i$, now for $i=1,2$. An important difference now is that the colors refer to the attractors. The units are distinguished by markers: circles for unit $i=1$ and diamonds for unit $i=2$. These markers correspond to the positions of the units at some arbitrarily chosen time point.

To begin this analysis, we recall that each uncoupled unit has three equilibria. Consequently, a system of two coupled units, under sufficiently weak coupling, has $3^2=9$ equilibria, corresponding to all combinations of the individual equilibria. Naturally, the symmetric combinations node-node, saddle-saddle, and focus-focus correspond to the two units being together in the same equilibrium. Since the coupling term becomes zero when the units are completely synchronized, these symmetric equilibria occupy the same positions as their uncoupled counterparts when projected into the units' subspace $x_i-y_i$. The other equilibria are asymmetric and have non-zero coupling terms, which shift their positions as a function of the coupling strength $\epsilon$. However, for simplicity, we still label the equilibria as combinations of the uncoupled equilibria, e.g. node-saddle denoting an equilibrium with $3$ negative eigenvalues and $1$ positive eigenvalue.

For $\epsilon = 0.05$, the node-node is the only attractor in the system (Fig. \ref{fig:attractors-2-units}A). In this solution, both units are in a steady state (SS), so we label the attractor as SS-SS (also called homogeneous steady state HSS \cite{ullner2008multistability}).

Next, at $\epsilon \approx 0.065$, a stable oscillation emerges, in which both units oscillate with a large amplitude (Fig.~\ref{fig:attractors-2-units}B). Therefore, we label this attractor LA-LA. It initially forms near the saddle point $\mathbf{x}_s^\mathrm{unc}$, located near the lower left corner. This proximity to the saddle point causes trajectories in that region to slow down significantly. As the coupling increases, the limit cycle moves farther away from the saddle point, resulting in a decreasing amplitude. This progression can be observed by comparing the attractors in subsequent panels.

At $\epsilon \approx 0.117485$ a pair of asymmetric attractors emerges, in which one unit has a large amplitude oscillation (LA) and the other unit has a small amplitude oscillation (SA), and vice-versa (Fig.~\ref{fig:attractors-2-units}C). Because the units are identical, the system has a permutation symmetry, so both attractors, LA-SA (large amplitude in unit 1 and small amplitude in unit 2), and SA-LA (reciprocal case) are simply permuted versions of each other. Consequently, these attractors overlap each other in Figs.~\ref{fig:attractors-2-units}C-D. They can be distinguished by the position of the units, indicated by the markers. Please note that the small amplitude oscillation has such a small amplitude that it is barely visible in the figures. In the literature, the LA-SA attractors have also been called inhomogeneous limit cycles \cite{ullner2008multistability, tyson1975control} (IHLC). 

At this coupling strength $\epsilon \approx 0.117485$, the system has four coexisting attractors, three of them being oscillations, even though the uncoupled dynamics only has equilibria! Eventually, for stronger coupling the pair LA-SA and SA-LA disappears around $\epsilon \approx 0.22$, and the system becomes bistable again. The result is shown in Fig.~\ref{fig:attractors-2-units}E. At $\epsilon \approx 0.27$, the periodic LA-LA attractor is replaced by a quasi-periodic LA-LA attractor, which again has both units oscillating with a large amplitude. In the quasi-periodic attractor, the units have different frequencies, and are desynchronized in both frequency and phase (cf. Fig.~\ref{fig:attractors-2-units}F). Eventually it disappears and only the stable equilibrium remains for sufficiently strong coupling. 
\begin{figure*}[ht!]
    \centering
    \includegraphics[width=\textwidth]{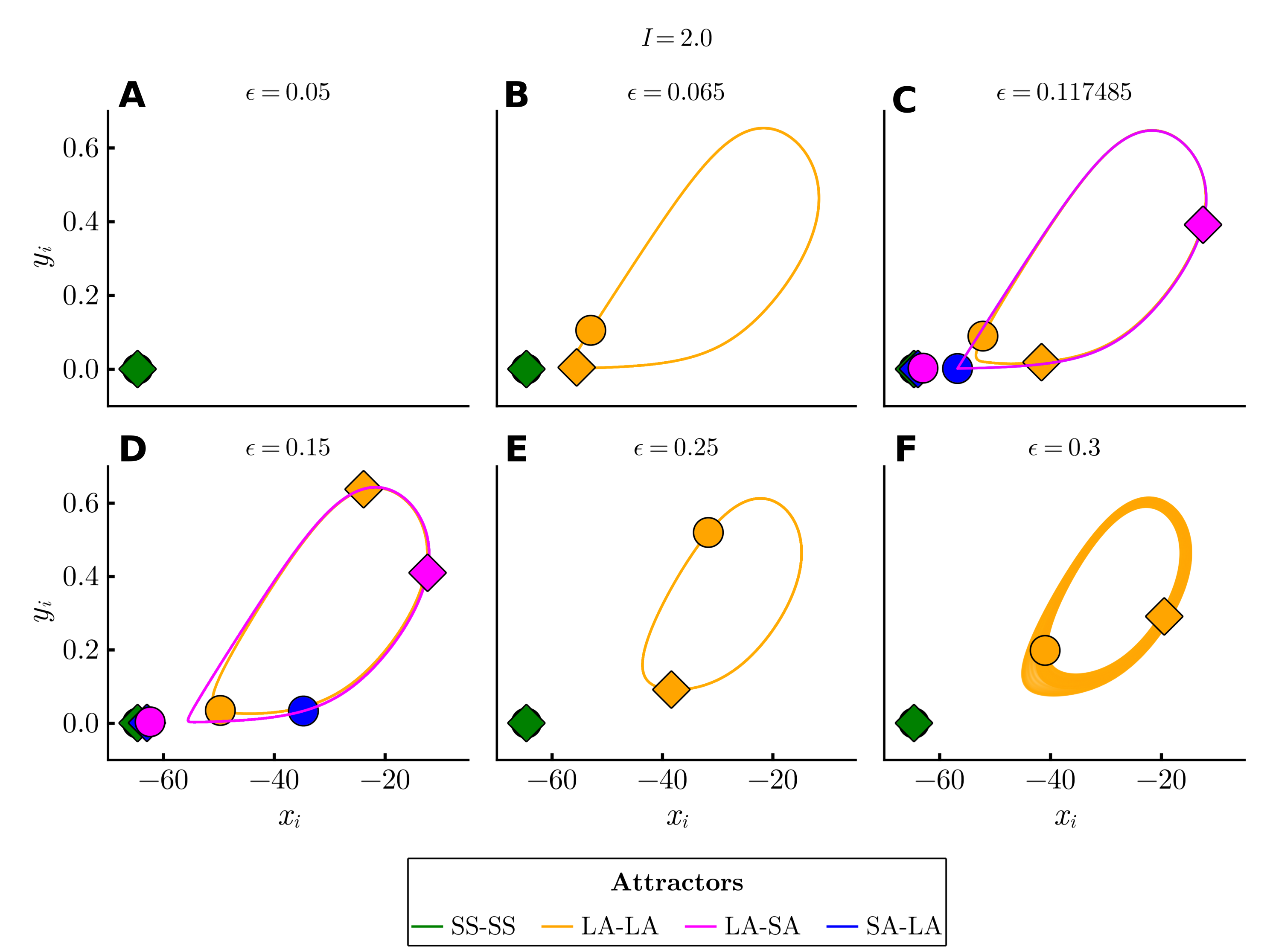}
    \caption{Attractors created by diffusive coupling of two coupled excitable units. Each panel is a projection onto $2D$ space of $x_i-y_i$ for different coupling strengths. Each attractor receives a unique color, defined in the legend. The markers denote the positions of the units for an arbitrarily chosen time point, with unit $i=1$ shown as a circle and unit $i=2$ as a diamond. The stable equilibrium SS-SS (green) is the only attractor existing for weak coupling strengths, as shown for $\epsilon=0.05$. Another attractor emerges at $\epsilon \sim 0.065$ corresponding to two units oscillating with large amplitude - it is thus labeled as LA-LA and colored yellow. A pair of asymmetric attractors emerges at $\epsilon \sim 0.117485$ corresponding to one unit oscillating with large amplitude and the other oscillating with small amplitude; they are labeled respectively as LA-SA (purple) and SA-LA (blue). The pair eventually disappears and the system becomes bistable again at $\epsilon=0.25$. At $\epsilon=0.3$, the LA-LA attractor is quasi-periodic. For stronger coupling $\epsilon$, the torus disappears, such that only the stable equilibrium is left for sufficiently strong $\epsilon$.}
    \label{fig:attractors-2-units}
\end{figure*}

\subsection{Bifurcations giving rise to the attractors}\label{sec:bifurcations}
To understand the emergence and disappearance of the periodic and quasiperiodic attractors in the $N=2$ case, we start by studying their associated bifurcations. We perform a continuation analysis using the XPPAUT 8.0 software \cite{ermentrout2002simulating}. This analysis is shown in Fig.~\ref{fig:bifurcations-2-units}, where the period $T$ of oscillation is estimated as a function of the coupling strength $\epsilon$. In this figure, the  green and red colors indicate stable and unstable solutions, respectively. First, in Fig.~\ref{fig:bifurcations-2-units}A, we present the continuation analysis for the LA-LA attractor, where both units oscillate with a large amplitude. We observe that this attractor arises from a saddle-node bifurcation of limit cycles (SNLC) at $\epsilon \approx 0.06432$. Subsequently, the stable limit cycle undergoes a Neimark-Sacker (torus) bifurcation (TR) at $\epsilon \approx 0.2701$, becoming unstable and being replaced by a stable torus. Next, this unstable limit cycle disappears in a supercritical Hopf bifurcation (HB) at $\epsilon \approx 0.4088$. Meanwhile, the saddle limit cycle that emerges at the SNLC bifurcation disappears in a homoclinic bifurcation (HOM) involving a saddle-saddle equilibrium at $\epsilon \approx 0.07285$. While it exists, the saddle limit cycle forms the basin boundary between the stable equilibrium and the stable limit cycle. When it disappears in the homoclinic bifurcation, it is immediately replaced by a pair of asymmetric saddle limit cycles that also emerges in a homoclinic bifurcation to the same equilibrium at the same parameter value, as shown in Fig.~\ref{fig:bifurcations-2-units}B. These saddle limit cycles then compose the basin boundary between the attractors. They correspond to the unstable version of the LA-SA and SA-LA attractors, which are later born also in a homoclinic bifurcation at $\epsilon \approx 0.1175$, but involving a saddle-node equilibrium. Eventually, both the stable and the unstable limit cycles collide and disappear in a SNLC bifurcation at $\epsilon \approx 0.2179$. The files used to perform the analysis are freely available at \onlinecite{rossi2024multistabilitygithub}.
\begin{figure}[h!]
    \centering
    \includegraphics[width=1.0\columnwidth]{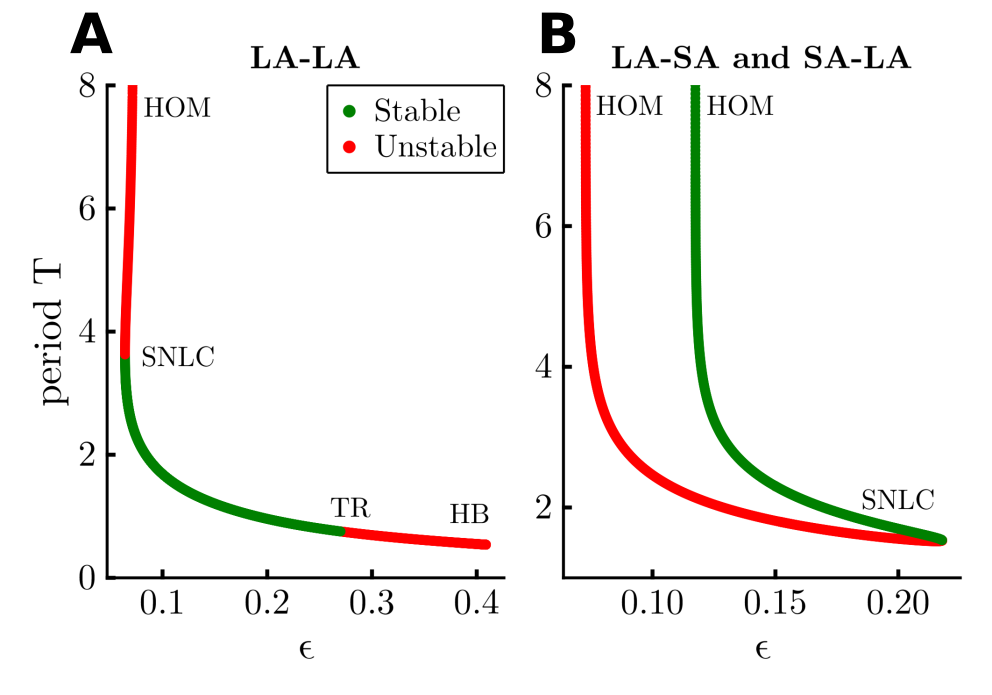}
    \caption{Continuation analysis for oscillations in two-unit case. Each panel shows a continuation of the limit cycles, plotting their period $T$ as a function of the coupling strength $\epsilon$ for fixed $I = 2.0$. The left panel shows the analysis for the LA-LA attractor, which emerges in a saddle-node bifurcation of limit cycles (SNLC) together with an saddle limit cycle (red curve) at $\epsilon \approx 0.06432$. The unstable LC goes through a homoclinic bifurcation (HOM) where it collides with a saddle-saddle equilibrium at $\epsilon \approx 0.07285$ and disappears. The LA-LA attractor (in green) remains stable until it loses stability due to a Neimarck-Sacker bifurcation (TR) and then disappears due to a supercritical Hopf bifurcation (HB). In panel B, the LA-SA, shown in the green curve, emerges due to a homoclinic bifurcation involving a saddle-node equilibrium at $\epsilon \approx 0.1175$ and then disappears due to a SNLC, when it collides with its unstable counterpart, in red, that is also born in a homoclinic bifurcation to a saddle-saddle equilibrium at $\epsilon \approx 0.07285$. The bifurcation diagram is identical for the SA-LA attractor, due to symmetry.}
    \label{fig:bifurcations-2-units}
\end{figure}

\subsection{Emergence of oscillations through reinjection mechanism}\label{sec:geometry}
To gain insights into how the coupling between the units supports the emergence of oscillations in a system whose uncoupled units exhibit only steady states, we now examine the geometry of the emerging attractors. As we see in Eqs.~\ref{eq:generic-network-equation}, the dynamics of unit $i$ can be decomposed into two terms: the local dynamics, governed by $\mathbf{f}_i(x_i,y_i)$, and the coupling, governed by $\epsilon \mathbf{h}_i(\mathbf{x})$. The local dynamics $\mathbf{f}_i(x_i, y_i)$ generates a vector field dictating the trajectories of the uncoupled units. As described in Sec.~\ref{sec:uncoupled}, $\mathbf{f}_i$ creates an excitability region, on which trajectories go through long excursions in state space before converging to the stable equilibrium. They follow the stable manifold $\mathbf{W}^s(\mathbf{x}_s^\mathrm{unc})$ of the saddle point $\mathbf{x}_s^\mathrm{unc}$ in the uncoupled system on their way to the equilibrium. The coupling dynamics $\mathbf{h}_i(\mathbf{x}) = (x_j-x_i, y_j-y_i)$ generates a vector field that points from unit $i$ to unit $j$, with an amplitude proportional to their distance. For $\epsilon > 0$, the coupling $\epsilon \mathbf{h}_i$ is attractive, as it pulls unit $i$ towards unit $j$. In the following examples, we see how interaction between these two terms leads to the emergence of the stable oscillations. 

\begin{figure*}[h!]
    \centering
    \includegraphics[width=0.95\textwidth]{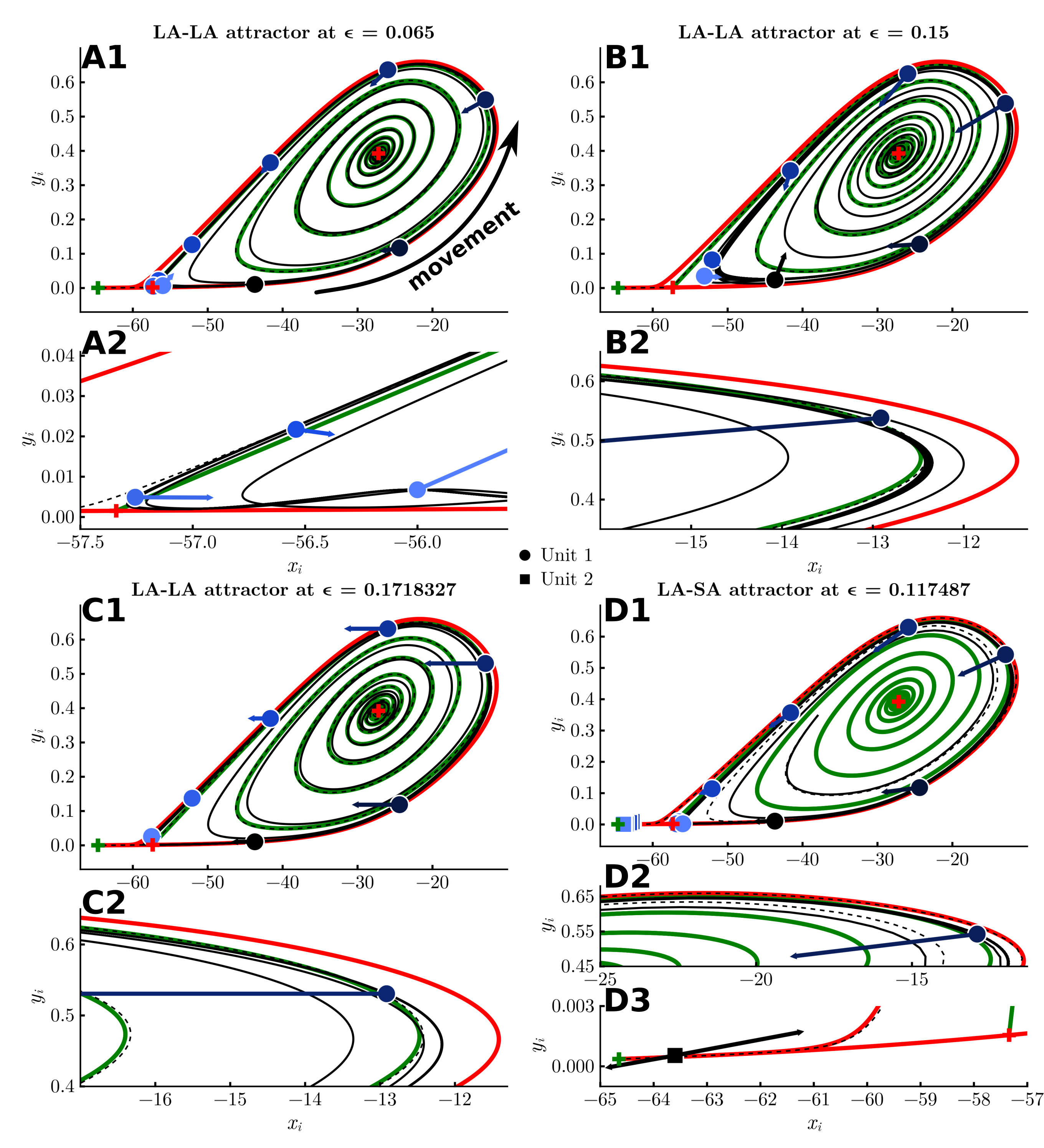}
    \caption{Illustration of the trapping phenomenon. Each panel shows a projection of the full 4D state space into the subspace $x_i$-$y_i$ of unit i. The three symmetric equilibria in the coupled system are shown: the stable equilibrium as a green cross, the unstable focus and the saddle point $\mathbf{x}_s^\mathrm{unc}$ as red crosses. The stable $\mathbf{W}^s(\mathbf{x}_s^\mathrm{unc})$ and unstable $\mathbf{W}^u(\mathbf{x}_s^\mathrm{unc})$ manifolds of the saddle $\mathbf{x}_s^\mathrm{unc}$ in the uncoupled system are also shown as green and red lines, respectively. The black solid lines represent an illustrative trajectory converging to one of the emerging attractors. A trajectory starting with the same initial condition but with $\epsilon=0$ is shown in black dashed lines. The position and coupling vector $\epsilon \mathbf{h}_1$ of unit 1 at specific time points are plotted respectively as circles and arrows. The unit's colors vary from black to light blue to indicate the passage of time. The LA-LA attractor is shown in panels A1-A2 for $\epsilon = 0.065$ and in B1-B2 for $\epsilon=0.15$. Panels C1-C2 show that this attractor still emerges if $\epsilon_2=0$, at the cost of requiring a larger value of $\epsilon_1$. Panels D1-D3 show LA-SA attractor, with the position of unit 2 (projected into $x_2-y_2$) also shown as squares. Panel D3 is added specifically to indicate the behavior of unit 2, with two arrows: the upward and rightward pointing arrow is the coupling $\epsilon \mathbf{h}_2$, and the leftward and downward arrow is the unit's uncoupled dynamics $\mathbf{f}_2$.}
    \label{fig:illustration_trapping}
\end{figure*}

Figure~\ref{fig:illustration_trapping}A1 illustrates this scenario for the LA-LA attractor at $\epsilon = 0.065$, already introduced in Fig.~\ref{fig:attractors-2-units}B. In Fig.~\ref{fig:illustration_trapping}A1, the structures in the complete 4D space are projected onto the $x_i-y_i$ subspace of each unit. For reference, we overlay on top of this plot the structures of the uncoupled unit, as seen already in Fig.~\ref{fig:excitable-neuron}. The stable equilibrium is represented as a green cross, while the saddle point $\mathbf{x}_s^\mathrm{unc}$ and the focus are shown as red crosses. The stable manifold $\mathbf{W}^s(\mathbf{x}_s^\mathrm{unc})$ of the saddle is a green line, while its unstable manifold $\mathbf{W}^u(\mathbf{x}_s^\mathrm{unc})$ is a red line. A trajectory converging to the LA-LA attractor is shown as a solid black line. Starting from an initial condition near the focus, at the center of the figure, the trajectory spirals outwards. This spiraling can be understood as the coupling being weak enough that the local dynamics $\mathbf{f}_i$ dominates the trajectory here, such that it roughly follows $\mathbf{W}^s(\mathbf{x}_s^\mathrm{unc})$.

Looking at an amplification of the region near $\mathbf{x}_s^\mathrm{unc}$ in Fig.~\ref{fig:illustration_trapping}A2, we see that the trajectory follows $\mathbf{W}^s(\mathbf{x}_s^\mathrm{unc})$ almost until the saddle $\mathbf{x}_s^\mathrm{unc}$ (red cross). Then, we see the crucial effect of the coupling. Without it, the trajectory would have followed along the left branch of $\mathbf{W}^u(\mathbf{x}_s^\mathrm{unc})$ and converged to the stable equilibrium. This is shown in the black dashed line, which shows a trajectory of the uncoupled system starting at the same initial condition (in $(x_1, y_1)$) as the black solid line. However, the coupled trajectory does not do that. Instead, it goes rightward, influenced by the coupling $\epsilon \mathbf{h}_i$. This effect can be seen by the rightward pointing arrows attached to the circles. The arrows correspond to the coupling vector $\epsilon \mathbf{h}_i$ on unit $i$, depicted as the circles whose colors vary along the trajectory, from dark blue to light blue. Under this coupling, the trajectory crosses $\mathbf{W}^s(\mathbf{x}_s^\mathrm{unc})$ in this projection, and is effectively reinjected into the excitability region.  For clarity, the trajectory of the coupled 4D system, when projected into the $x_i-y_i$ plane, crosses the stable manifold $\mathbf{W}^s(\mathbf{x}_s^\mathrm{unc})$ of the saddle of the uncoupled $2D$ system. Naturally, it does not cross any invariant manifolds of the coupled system.

As this attractor is symmetric, the behavior described for unit $1$ occurs identically for unit $2$. The reinjection into the excitability region thus happens for both units, causing them to repeatedly pull on each other and reinject each other into the previously transient region. In this sense, we say that \textit{the coupling traps the units in the excitability region}, preventing them from following their local dynamics' tendency toward the stable equilibrium.

Why does this crossing happen so close to $\mathbf{x}_s^\mathrm{unc}$? As we have seen in Sec.~\ref{sec:bifurcations}, the LA-LA attractor emerges in a saddle-node bifurcation of limit cycles (SNLC), and is not directly related to $\mathbf{x}_s^\mathrm{unc}$. However, the local dynamics has a magnitude $|\mathbf{f}_i|$ that is small in the vicinity of $\mathbf{x}_s^\mathrm{unc}$. So in this region the relative effect of the coupling $\epsilon \mathbf{h}_i$ increases. Near $\mathbf{x}_s^\mathrm{unc}$, the unit tends to move very slowly due to its local dynamics, but at the same time the coupling intensity is relatively strong. As a result, the coupling overcomes the local dynamics and unit $i$'s trajectory moves to the right near $\mathbf{x}_s^\mathrm{unc}$ on the $x_i-y_i$ projection. In summary, the slowness near $\mathbf{x}_s^\mathrm{unc}$ helps the coupling $\epsilon \mathbf{h}_i$ to overcome the local dynamics. 

The slowness near $\mathbf{x}_s^\mathrm{unc}$ also allows us to understand the bulge that occurs right after the trajectory crosses $\mathbf{W}^s(\mathbf{x}_s^\mathrm{unc})$, on the right-hand side of Fig.~\ref{fig:illustration_trapping}A2. The trajectory in this region is quite slow. In fact, unit $i$ spends most of its time on it while unit $j$ traverses the rest of the oscillation. As $j$ moves around the oscillation, the distance between the units increases significantly, and therefore so does $\mathbf{h}_i$ (note the longer blue arrow for unit 1 in the bulge). The combination of the slowness of $\mathbf{f}_i$ and the high value of $\mathbf{h}_i$ means that the coupling dominates the sum, significantly impacting the trajectory, pulling it towards unit $j$, creating the upwards movement of the bulge for unit $i$ (and vice-versa for unit $j$, because of the symmetry). 

For bigger values of $\epsilon$, the coupling becomes stronger and the slowness near $\mathbf{x}_s^\mathrm{unc}$ becomes less relevant. In Fig.~\ref{fig:illustration_trapping}B1, note how the coupling is larger (longer arrows) for $\epsilon = 0.15$. Consequently, the coupling manages to pull unit $i$, in this $x_1-y_1$ projection, across $\mathbf{W}^s(\mathbf{x}_s^\mathrm{unc})$ earlier along the manifold - see the magnification in Fig.~\ref{fig:illustration_trapping}B2. 
The stronger coupling also affects the shape of the attractor. This is most visible close to $\mathbf{x}_s^\mathrm{unc}$, which is considerably shifted rightward and upward if compared to the smaller $\epsilon$ in Fig.\ref{fig:illustration_trapping}A1-A2. This region has the slowest dynamics, and is thus the one most sensitive to the coupling. Furthermore, similarly to the argument leading to the bulge in Fig.\ref{fig:illustration_trapping}A1-A2, while unit $i$ is in this slow region, unit $j$ eventually becomes diametrically opposite it, and the coupling amplitude grows significantly. This increase, combined with the slow dynamics, pulls the units upward and rightward, explaining the shift. This decreases the amplitude of the oscillation, consistent with the behavior seen for larger networks (Figs.~\ref{fig:multistability-network}B-H), in which the amplitude of the oscillation is inversely proportional to the number of neighbors a unit has.

So far we have considered what happens in the case that the $x$ and $y$ directions are coupled with the same intensity, i.e., when $\epsilon = \epsilon_1 = \epsilon_2$. However, this is not required for new attractors to emerge. In particular, the LA-LA attractor still emerges when only the $x$-component of the coupling is kept (i.e., when $\epsilon_1=\epsilon$ and $\epsilon_2 = 0$). The reinjection occurs similarly to before, as shown in Figs.~\ref{fig:illustration_trapping}C1-C2, where an illustrative trajectory of the coupled system can again be seen to cross $\mathbf{W}^s(\mathbf{x}_s^\mathrm{unc})$. In fact, the $y$-component of the coupling is not necessary to generate the LA-LA attractor, although it helps. Decreasing $\epsilon_2$ from $\epsilon_2=\epsilon$ to $\epsilon_2=0$ has the effect of increasing the critical value of $\epsilon_1$ for the saddle-node bifurcation of limit cycles that creates the attractor, effectively postponing its emergence, but not inhibiting it. The example in Figs.~\ref{fig:illustration_trapping}C1-C2 occurs soon after the LA-LA attractor emerges. Note that the coupling is much bigger than it was for Figs.~\ref{fig:illustration_trapping}A1-A2. Conversely, decreasing $\epsilon_1$ from $\epsilon_1 = \epsilon$ to $\epsilon_1 = 0$ can either destroy the LA-LA attractor through a SNLC bifurcation or cause it to lose stability. Therefore, the $x$-component is necessary to generate this attractor. 

In Fig.~\ref{fig:illustration_trapping}D1 we return to the case $\epsilon_1=\epsilon_2=\epsilon$ and examine the asymmetric attractor LA-SA. Its geometry differs from the previous cases, since now there is an asymmetry between the units. As in the previous cases, unit 1 oscillates in the excitability region with a large amplitude. Meanwhile, unit 2 (squares) is positioned, in the $x_2-y_2$ projection, between the stable equilibrium (green cross) and the saddle point $\mathbf{x}_s^\mathrm{unc}$ of the uncoupled dynamics (red cross). Both units can be seen in Figs.\ref{fig:attractors-2-units}D1, represented respectively as circles and squares. Their colors denote different time points along the trajectory. In this configuration, unit $1$ is pulled downwards and to the left by unit 2, as illustrated in Fig.\ref{fig:attractors-2-units}D2, which shows $\epsilon \mathbf{h}_1$. This pull is capable of causing unit $1$ to cross $\mathbf{W}^s(\mathbf{x}_s^\mathrm{unc})$ and to be reinjected into the excitability region, as shown in Figs.~\ref{fig:illustration_trapping}D1-D2. Meanwhile, unit 2 is pulled rightwards and upwards by unit 1. This pull is counteracted by the attraction it feels towards the stable equilibrium, with the result being a small-amplitude oscillation. This competition is illustrated in Fig.\ref{fig:attractors-2-units}D3, where two arrows are associated with unit 2 for a representative time point of the trajectory: the upward and rightward arrow is $\epsilon \mathbf{h}_2$, which pulls unit 2 towards unit 1; the downward and leftward arrow is $\mathbf{f_2}$, which pulls unit 2 towards the stable equilibrium. 

With this configuration of the units in the LA-SA attractor, the $x$-direction is actually counter-productive. To see this, we can focus on the region to the right of the saddle point - in Fig.\ref{fig:attractors-2-units}D1, several instances of unit 1 can be seen accumulated in this region with light blue colors. This accumulation is a result of the slow dynamics of trajectories moving near the saddle point. In this region, unit 1 is pulled leftwards back towards $\mathbf{W}^s(\mathbf{x}_s^\mathrm{unc})$. The $x$-direction is thus acting against the reinjection mechanism, and so may be impeding the emergence of the LA-SA attractor. We can verify this qualitative claim by decreasing $\epsilon_1$ from $\epsilon_1=\epsilon$ towards $\epsilon_1 = 0$. By doing this, we see that indeed the critical value of $\epsilon_2$ that leads to the emergence of the attractor decreases. If $\epsilon_1=0$ the attractor can still emerge. Therefore, the LA-LA and LA-SA attractors exhibit different dependencies on the $x$- and $y$-components of the coupling, due to their distinct geometries. This distinction is crucial in various applications where coupling can model diverse phenomena. For example, in ecology, coupling in the $x$-direction might represent migration between prey species, while coupling in the $y$-direction could denote migration among predator species.

Another intriguing attractor is the torus that emerges from the LA-LA attractor (cf. Fig.~\ref{fig:attractors-2-units}F), in which the two units continue to oscillate with a large amplitude, but quasi-periodically. Its geometry resembles that of the LA-LA attractor, but the increased coupling strength causes the units to exert a stronger mutual influence. Consequently, the reinjection still occurs, but now one unit pulls on the other so forcefully that the trajectory no longer follows a closed curve. Indeed, as $\epsilon$ increases, the torus expands, becoming wider.

Therefore, the reinjection mechanism, which acts to trap units in the excitability region of their uncoupled dynamics, underlies the emerging attractors we observe. This is true for the two distinct geometries: LA-LA and LA-SA, which emerge from different bifurcations. The mechanism also occurs for different dynamics: periodic and quasi-periodic.

\section{Discussion}\label{sec:discussions}
In this paper we have shown that diffusive coupling acting on excitable dynamics can create multistability of oscillations. The variety of coexisting attractors, which can be periodic, quasiperiodic, and even chaotic, emerge in a similar way: through the trapping of units in the excitability region of their local dynamics. This local dynamics consists of three equilibria living in the units' state space: an unstable focus, a saddle point $\mathbf{x}_s^\mathrm{unc}$ and a stable equilibrium, which is the only attractor in the uncoupled system. The stable manifold $\mathbf{W}^s(\mathbf{x}_s^\mathrm{unc})$ of $\mathbf{x}_s^\mathrm{unc}$ is extended in state space and has one branch that spirals out of the unstable focus. As a consequence, it separates the nearby state space into two regions: one that directly converges to the stable node and another that has to go on a long excursion around the stable manifold before converging to the node. On top of this, the units feel the attractive diffusive coupling, which pulls one unit toward the other. The dynamics of the coupled units is determined by the interaction of these two effects: the local dynamics attempting to pull the units towards the stable node and the diffusive coupling attempting to pull the units towards each other.

This competition is controlled by the coupling strength $\epsilon$. As already described in the literature for similar systems, there are two extremes \cite{stankovski2017coupling}. For sufficiently small $\epsilon$, the local dynamics dominates, and the only attractor is the stable node. For sufficiently large $\epsilon$ the coupling dominates, and the units converge to each other. When they do so, the coupling becomes zero, and then they again follow their uncoupled dynamics and converge to the stable node only. It is in between these extremes that interesting dynamics can occur \cite{stankovski2017coupling}. In this case, the coupling is strong enough to impact the trajectory of the uncoupled system, but not enough to completely overrule it. Because of the geometry of state space, the coupling can manage to pull the units away from the stable node and into the excitability region. The units find a stable configuration in which they are repeatedly reinjected into the excitability region, generating permanent oscillations. The type of these oscillations depends on the coupling strength, the number of interacting units, and the network's topology.

It has been known that diffusive coupling on units with a single stable equilibrium in a region of state space can create oscillations. These oscillations can be periodic, originating from a Hopf bifurcation of the equilibrium \cite{Smale1976, pogromsky1999on} or chaotic oscillations \cite{kocarev1995on} originating from a Shilnikov homoclinic bifurcation \cite{nijholt2023chaotic}. However, as we have shown, the scenario for an excitable system, with the additional interaction of two unstable equilibria, has important differences to the single equilibrium case. First, for $N=2$ coupled units, we have shown that periodic attractors can coexist with other periodic attractors and with the stable equilibrium, leading to a multistable coupled system. Further, these periodic attractors are qualitatively different: in one attractor, both units oscillate with a large amplitude (LA-LA attractor); in the other attractor, one unit instead has a very small amplitude oscillation, almost stationary (LA-SA and SA-LA attractors). The $N=2$ case also supports the emergence of a quasiperiodic oscillation, which coexists with the stable equilibrium. The bifurcations giving rise to the periodic attractors also differ: the periodic attractors emerge either through a saddle-node bifurcation of limit cycles (SNLC) bifurcation or through a homoclinic (HOM) bifurcation.

In the bigger network, with $N=10$ units, we have shown that the multistability becomes even richer. In this case, all these types of dynamics can coexist. The sheer number of coexisting attractors is also large (we observed up to 84 attractors for only $N=10$ units), with a dominance of periodic solutions. Since units can be either trapped in the excitability region, with large amplitude oscillations, or oscillate with low amplitude near $\mathbf{x}_s^\mathrm{unc}$, adding more units leads to a higher number of possible combinations of which units are placed in which position. Not all of these combinations are necessarily invariant solutions; and the ones that are invariant are not necessarily stable. The invariance and stability are controlled by the topology of the network, and more research in the future is needed to understand how exactly. It would be interesting to understand which topologies maximize the number of attractors, and which minimize them. This could provide further insights into other systems, for which a scaling of the number of attractors with the size of the network has been observed \cite{ullner2007multistability, rossi2022shifts, gelbrecht2020monte}.

By studying the attractors from the point of view of the local dynamics competing with the coupling, we identified an impact that the topology has on the attractors of the coupled system. Units that receive more connections tend to have a stronger coupling term than units with fewer connections. A stronger coupling coupling term pulls the units more strongly towards the excitability region. This effect is stronger in the slower region of the oscillation, close to $\mathbf{x}_s^\mathrm{unc}$, where the local dynamics is weaker. Thus, one could expect the oscillation of the units with more connections to be pushed away further from $\mathbf{x}_s^\mathrm{unc}$, and thus have a smaller amplitude than units with fewer connections. This is indeed what we observe: units with more neighbors have smaller amplitudes.

On the attractors, the units are permanently reinjected into the excitability region by the coupling. In this sense, they are trapped in a transient region (transient for the uncoupled dynamics). Trapping in transient regions due to coupling appears to be a common mechanism for creating new attractors in networked systems. We have observed a similar behavior in an excitable model of an ecological predator-prey system based on the Truscott-Brindley model \cite{truscott1984ocean}. There, new equilibria are created by the coupling. On the equilibria, the coupling exactly balances out the local dynamics, and the units reach an equilibrium. Another example of trapping has been elucidated in units with chaotic saddles in their local dynamics \cite{medeiros2018boundaries, medeiros2019state, medeiros2021the}. Chaotic saddles are non-attracting invariant chaotic sets. Under some circumstances, when a sufficiently large number $N$ of units are coupled diffusively, they can get trapped in this chaotic saddle and form a chaotic attractor \cite{medeiros2018boundaries, medeiros2019state}.  Another example has been recently elucidated in a system with a canard \cite{contreras2023scale}. 

Interestingly, in the neuronal system we have studied, the trapping also works if the coupling is present in only one of the directions $x$ or $y$. These directions have different effects in generating new attractors, due to the geometry of state space. If only the $x$-direction of the coupling is present, only the LA-LA attractor emerges, not the LA-SA or SA-LA. For coupling in the $y$-direction the reverse is true: LA-LA does not emerge, but LA-SA and SA-LA do. This coupling in $y$ is not biophysically relevant for the neuronal system we study here, but is important in ecological systems. There, this may represent the difference between migration of predators or the migration of prey species. In particular for the $x$-direction, we also mention that multistability is still maintained in a bigger networks with $N=10$ units. The coexisting oscillating attractors, with a mixture of some units oscillating at large amplitude and some at low amplitudes still occur, with more units in high-amplitude oscillations than in the case with the $x$ and $y$-couplings. We believe our study could serve as inspiration for future studies in other systems, such as in ecological ones, to investigate these effects of the coupling in more detail. Furthermore, it also serves as simple yet powerful example of the more general phenomenon of multistability through trapping of units in transients.


\section*{Acknowledgment}
K.L.R. was supported by the German Academic Exchange Service (DAAD).  E.S.M and U.F. acknowledge the support by the Deutsche Forschungsgemeinschaft (DFG) via the project number 454054251 (FE 359/22-1) and by The S\~ao Paulo Research Foundation (FAPESP) via the project number 2023/15040-0. The simulations were performed at the HPC Cluster ROSA, located at the Carl von Ossietzky Universität Oldenburg (Germany) and funded by the DFG through its Major Research Instrumentation Programme (INST 184/225-1 FUGG) and the Ministry of Science and Culture (MWK) of the state of Lower Saxony.

\appendix

\section{Supplemental material}

\subsection{Attractors for different local dynamics of the neurons}
In the main text we study the behavior of the coupled neuronal units with fixed parameters. One important parameter we can vary is $I$. Taking it as a bifurcation parameter and increasing it, the manifolds of the saddle approach each other to form a homoclinic orbit at $I = I_\mathrm{HOM} \approx 3.09$. At $I>I_\mathrm{HOM}$, a stable limit cycle emerges from the homoclinic orbit, and now the unit is bistable. Increasing $I$ further, the saddle and the node approach each other and a saddle-node bifurcation occurs at $I = I_\mathrm{SN} \approx 4.8$, and the neuron goes back to being monostable, with only the stable limit cycle remaining.

\begin{figure*}[h!]
    \centering
    \includegraphics[width=0.8\textwidth]{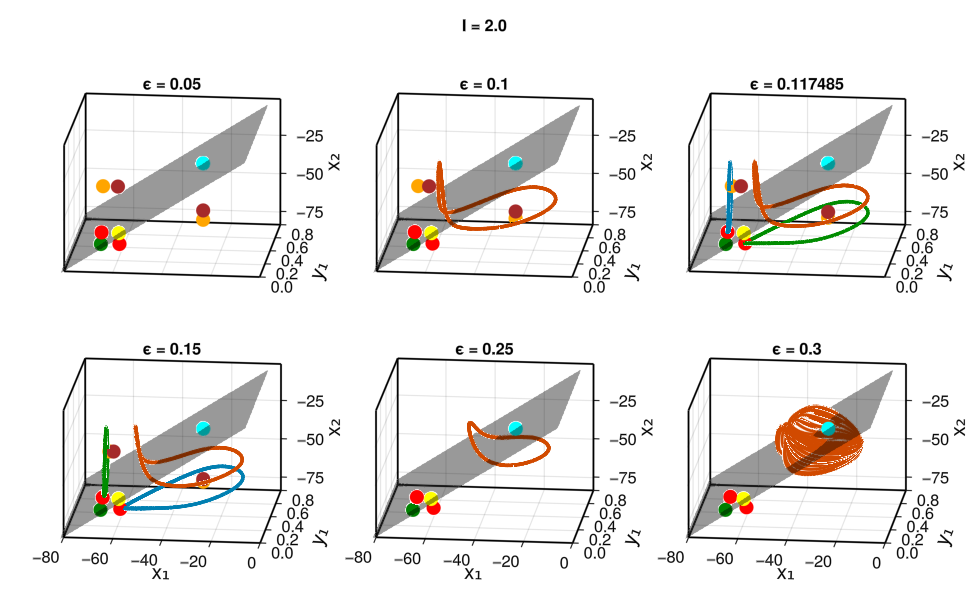}
    \caption{Attractors transformed or created by the diffusive coupling for $I=2.0$. Each panel is a projection onto $3D$ space spanned by $x1-y1-x2$. The gray plane denotes the plane with $x_1=x_2$. The circles' colors denote types of equilibria. These equilibria are formed as combinations of the three equilibria in the single units, and are labeled according to this combination: green for stable node, yellow for saddle-saddle (two positive, two negative eigenvalues), cyan for focus-focus (two pairs of complex conjugate eigenvalues with positive real part), red for saddle-node (three positive, one negative eigenvalue), orange for node-focus (one pair of complex conjugate eigenvalues with positive real part, two negative eigenvalues) and brown for saddle-focus (one pair of complex conjugate eigenvalues with positive real part, one positive and one negative eigenvalues). One attractor emerges at $\epsilon \sim 0.065$ corresponding to two units oscillating with large amplitude. Two (symmetric) attractors emerge at $\epsilon \sim 0.117485$ corresponding to one unit oscillating with small amplitude around the saddle-node point and the other oscillating with large amplitude. The symmetric attractors die out at $\epsilon \approx 0.22$. A torus emerges at $\epsilon \sim 0.25$. More bifurcations keep happening until all new attractors die out, and only the stable equilibrium is left. }
    \label{fig:attractors-2-units-3d}
\end{figure*}

To recall the attractors emerging at $I=2.0$ and to provide a complementary view, we show in Fig. \ref{fig:attractors-2-units-3d} a three-dimensional version of Fig. \ref{fig:attractors-2-units}. One can notice the emergence of the LA-LA attractor, and its eventual replacement by a torus. Also one can see the emergence of the LA-SA attractors, which emerge touching the red circles, that denote saddle points of saddle-node type. 

\begin{figure}[htb!]
    \centering
    \includegraphics[width=\columnwidth]{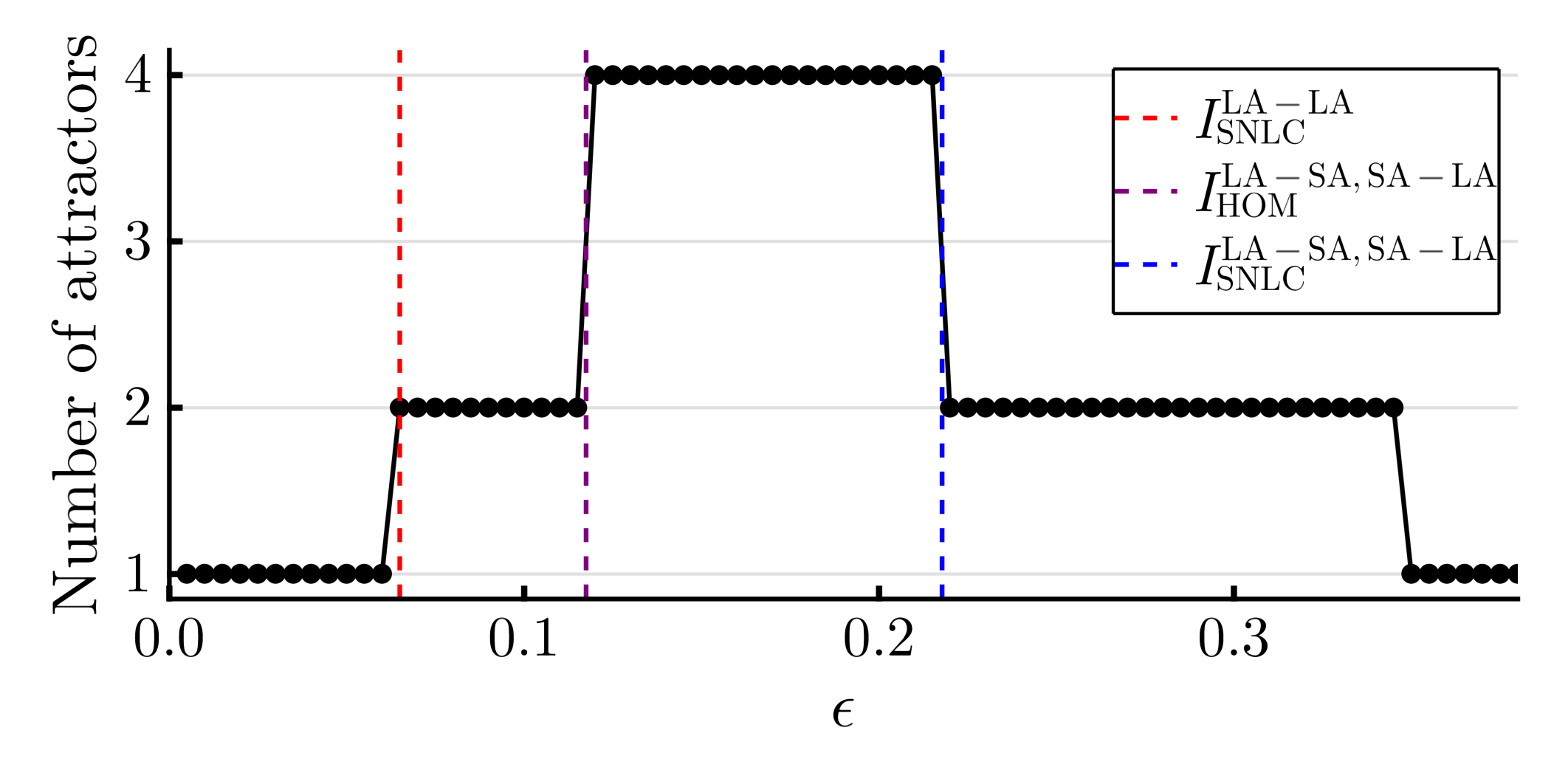}
    \caption{Number of attractors found for two excitable units at $I=2.0$. The red, purple, and blue dashed lines respectively represent the coupling strength values $\epsilon$ at the following bifurcations occurs: the saddle-node bifurcation of limit cycles (SNLC) generating the large amplitude-large amplitude (LA-LA) attractor; the homoclinic (HOM) bifurcation generating the large amplitude-small amplitude (LA-SA) attractor and its symmetric pair; and the SNLC bifurcation destroying this pair.}
    \label{fig:number-attractors}
\end{figure}

To summarize the multistability picture, we show the number of coexisting attractors we numerically identified for each value of coupling strength $\epsilon$ in Fig.~\ref{fig:number-attractors}. For sufficiently weak coupling strength $\epsilon$ only one attractor, the stable equilibrum, is present. Then a saddle-node bifurcation of limit cycles (SNLC) occurs that generates the large amplitude-large amplitude (LA-LA) attractor, marked by the red dashed lines, and the system becomes bistable. Eventually a homoclinic bifurcation (HOM) occurs (purple dashed line) generating the large amplitude-small amplitude (LA-SA) attractor and its symmetric reciprocal SA-LA. Now, the system has four coexisting attractors. Subsequently, the LA-SA and SA-LA attractors disappear in a SNLC bifurcation (blue dashed lines). Eventually the LA-LA attractor becomes unstable, replaced by a stable torus which later also disappears. For strong $\epsilon$ in the considered interval only the stable equilibrium is left.

\begin{figure}[h!]
    \centering
    \includegraphics[width=\columnwidth]{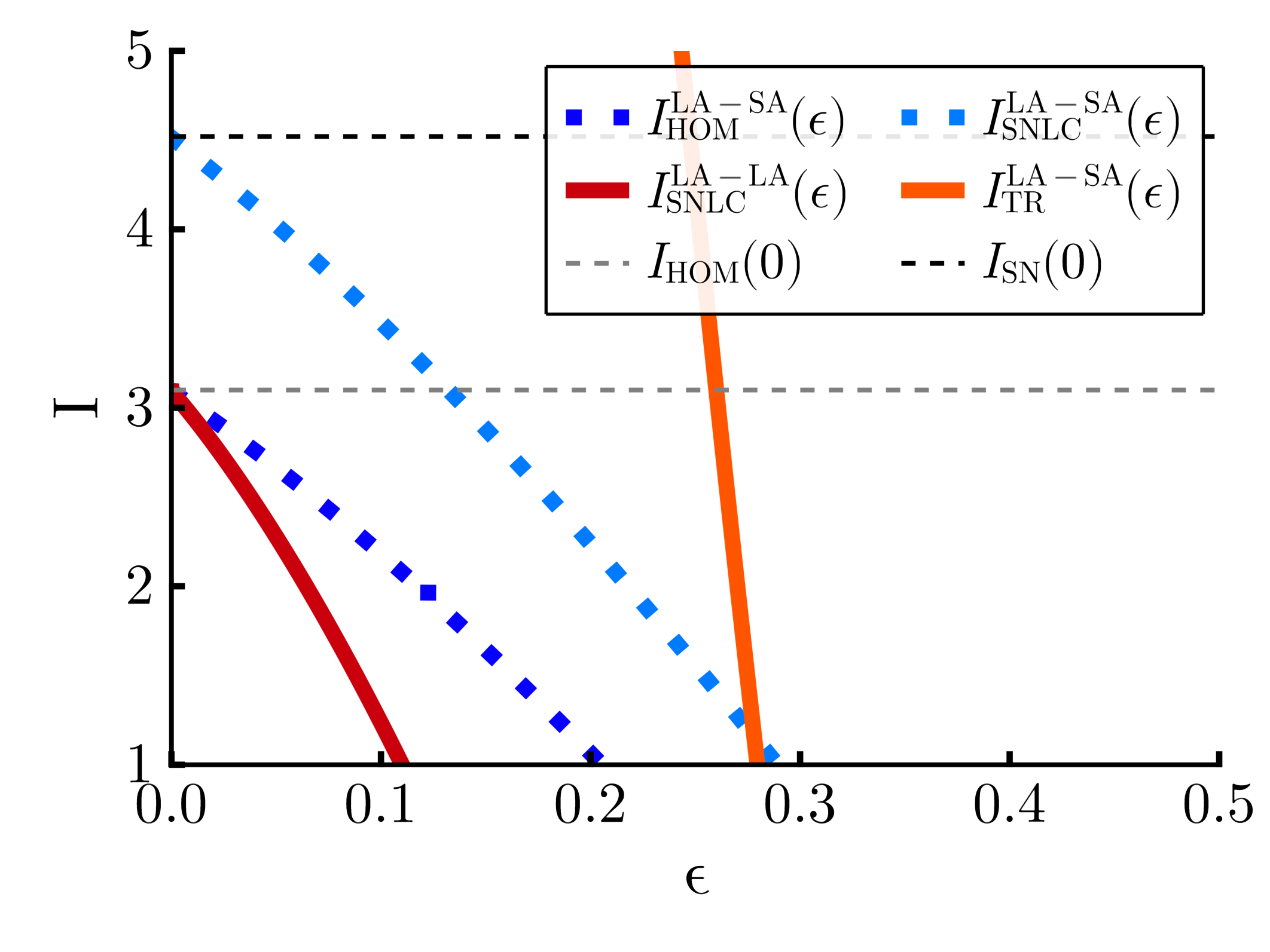}
    \caption{Two-parameter continuation curves across $I$ and $\epsilon$. The curves denote the $I(\epsilon)$ combinations that lead to each labeled bifurcation. The red and orange solid curves denote the bifurcations for LA-LA attractors, born through a SNLC bifurcation and de-stabilized through a torus (TR) bifurcation. The LA-LA is thus stable in between those curves. The blue and cyan dotted curves denote bifurcations occurring for the LA-SA attractor, born through a homoclinic (HOM) bifurcation and disappearing through a saddle-node bifurcation of limit cycles (SNLC). The LA-SA exists in between those curves. The homoclinic and saddle-node of equilibria (SN) bifurcations occurring in the uncoupled ($\epsilon=0$) case are respectively shown in grey and black dashed lines.  Continuations were done using XPPAUT \cite{ermentrout2002simulating}}.
    \label{fig:twoparam-continuation}
\end{figure}

To understand how these attractors depend on the local dynamics of the units, particularly how they change when the units go through the homoclinic and saddle-node bifurcation, we have studied a two-parameter continuation curve of the bifurcations giving rise to both the LA-LA and the LA-SA attractors. This is shown in Fig. \ref{fig:twoparam-continuation}.

The SNLC bifurcation generating the LA-LA solution converges to $(\epsilon, I) = (0, I_\mathrm{HOM})$. For $I>I_\mathrm{HOM}$ the LA-LA seems to occur for any value $\epsilon > 0$ that we tested. Therefore it seems that the LA-LA SNLC curve becomes vertical for $I>I_\mathrm{HOM}$ at $\epsilon = 0$.
The SNLC bifurcation destroying the LA-SA attractor converges to $(\epsilon, I) = (0, I_\mathrm{SN})$.

The attractors therefore exist for a wide range of parameters in the local dynamics of the units. Further, for $I>I_\mathrm{HOM}$, another attractor emerges, in which both units synchronize completely in the limit cycle that is now stable in the local dynamics.

\end{document}